\documentclass [12 pt]{article}
\usepackage {a4}
\usepackage{amsmath}
\usepackage {amssymb}

\begin{document}

\title{On the relation between lifting obstructions
and ordinary obstructions} 

\author{Christian Bohr \\ 
Mathematisches Institut \\ Theresienstr. 39 \\ 80333 M\"unchen \\
bohr@rz.mathematik.uni-muenchen.de
}
\date{\today}
\maketitle

\setlength{\unitlength}{0.1 cm}

\newsavebox{\diagrammeins}
\savebox{\diagrammeins}{
\put(0,20){\dots} \put(6,20){\vector(1,0){10}}
\put(17,19){$B^{k+1}$} \put(26,20){\vector(1,0){10}}
\put(28,21){$\scriptstyle f_{k+1}$}
\put(39,19){$B^k$} \put(46,20){\vector(1,0){10}}
\put(50,21){$\scriptstyle f_k$}
\put(57,19){$B^{k-1}$} \put(66,20){\vector(1,0){10}}
\put(77,20){\dots}
\put(39,0){$B$}
\put(38,3){\vector(-1,1){15}}
\put(44,3){\vector(1,1){15}}
\put(41,4){\vector(0,1){13}}
\put(38,8){$\scriptstyle \iota_k$}
\put(30,12){$\scriptstyle \iota_{k+1}$}
\put(46,12){$\scriptstyle \iota_{k-1}$}
}

\newsavebox{\diagrammzwei}
\savebox{\diagrammzwei}{
\put(0,22){$E$} \put (25,22){$P$}
\put(5,23){\vector(1,0){18}}
\put(0,2){$B$} \put(23,2){$B^{q-1}$}
\put(5,3){\vector(1,0){17}}
\put(2,20){\vector(0,-1){12}}
\put(26,20){\vector(0,-1){12}}
\put(3,14){$\scriptstyle p$}
\put(27,14){$\scriptstyle p'$}
\put(12,1){$\scriptstyle \iota_{q-1}$}
\put(14,24){$\scriptstyle \eta$}
\put(36.5,13){$B^q$}
\put(5,4){\line(3,1){19}}
\put(27,11){\vector(3,1){8}}
\put(35,12){\vector(-1,-1){6}}
\put(14,9){$\scriptstyle \iota_q$}
}

\newsavebox{\diagrammdrei}
\savebox{\diagrammdrei}{
\put(0,35){$E$}
\put(11,24){$S$}
\put(3,34){\vector(1,-1){7}}
\put(1,34){\vector(1,-3){9}}
\put(10.5,2){$B$}
\put(12.5,23){\vector(0,-1){17}}
\put(23,35){$\tilde{E}$}
\put(22,34){\vector(-1,-1){7}}
\put(24,34){\vector(-1,-3){9}}
\put(5,36){\vector(1,0){16}}
\put(12,37){$\scriptstyle \lambda$}
\put(13,15){$\scriptstyle \pi$}
\put(1.5,25){$\scriptstyle p$}
\put(22.5,25){$\scriptstyle \pi \circ q'$}
\put(7,31){$\scriptstyle q$}
\put(17,31){$\scriptstyle q'$}
}

\newtheorem{satz}{Proposition}
\newtheorem{theorem}{Theorem}
\newcommand {\gammaquer}{\overline{\gamma}}
\newcommand {\cquer}{\overline{c}}
\newcommand {\fquer}{\overline{f}}
\newcommand {\gquer}{\overline{g}}
\newcommand {\dquer}{\overline{d}}
\newcommand {\hdach}{\hat{h}}
\newcommand {\hquer}{\overline{h}}
\newcommand {\Z}{\mathbb{Z}}
\newcommand {\R}{\mathbb{R}}

\noindent
{\bf \small Abstract:} {\it \small
We consider partial liftings $k : A \rightarrow E$ of maps
$f : X \rightarrow B$, where $(X,A)$ is a relative CW-complex and $E
\rightarrow  B$ is a fibration. In this situation, we have a primary
obstruction to extend the partial lifting to a lifting of $f$ on all of $X$,
and there is an obstruction to extend $k$ as an ordinary map into the space
$E$. A relation between these two cohomology classes is proved when the 
fibre of
$E \rightarrow B$ is an Eilenberg-Mac Lane space $K(\Pi,n)$ and $\pi_i(E)=0$ 
for
$i \leq q-1$, where $q \geq n+2$, that specialises to well-known formulas
about secondary obstructions. The result is applied
to the Hopf fibration (what includes defect sections in $S^1$-bundles
over 4-manifolds)
and to the case of a certain $SU(3)$-bundle over $S^4$. }



\section {Introduction and summary of results}

\noindent
For a first example, let $M$ be a 4-dimensional closed and 
oriented manifold and 
$$
f : M \backslash \Delta \rightarrow S^2
$$
a continous map, where $\Delta \subset M$ is a finite
set. Following the terminology that is used in the study of defects
in ordered media, $f$ will be called a 
``defect map" and $\Delta$ the set of  ``point defects" 
(see \cite{M},\cite{J} 
and in particular \cite{KT} for the physical aspects of the questions
treated in this paper).
To each point $p \in \Delta$ we
can assign its {\bf local index} $\iota_p(f) \in \pi_3(S^2)=\mathbb{Z}$, given
by the restriction of $f$ to the boundary of an embedded ball $D^4 \subset M$
whose center is $p$ and whose intersection with $\Delta$ consists of $p$ alone
(the embedding is supposed to preserve the orientation).
Now consider the Hopf fibration $S^3 \rightarrow S^2$ and a defect lifting
$$
\overline{f} : M \backslash (\Delta \cup \Delta') \rightarrow S^3
$$
of $f$. The defect set $\Delta'$ can be chosen to be a closed orientable
surface
in $M \backslash \Delta$: Let $\xi:=f^*(S^3 \rightarrow S^2)$ denote the
induced $S^1$-bundle over $M \backslash \Delta$ and choose a transversal
section $\sigma : M \backslash \Delta \rightarrow \xi \times _{S^1} \mathbb
{C}$ in the associated complex line bundle that is not zero on a certain 
neighborhood of $\Delta$. A defect lifting with defect
set $\Delta'=\sigma^{-1}(0)$ is now given by $\frac{\sigma}{\| \sigma \|}$.

Again, we have local indices that reflect the behavior of
$\overline{f}$ near $\Delta'$: Let $\Delta'=\bigcup_i \Delta_i'$ be the
decomposition of $\Delta'$ into its connected components. For each $i$ choose
a point $\delta_i \in \Delta_i'$ and an embedding $D^2 \rightarrow M
\backslash (\Delta \cup \bigcup _{j \neq i} \Delta_j')$, given by the
restriction of a tubular 
neighborhood of $\Delta_i'$ to the fibre of the normal bundle of $\Delta_i'$
over $\delta_i$. Since the induced bundle $\xi$ is trivial over $D^2$, the 
lifting $\overline{f}$ (that is just a section of $\xi$) on the boundary
of $D^2$ is given by a map $S^1 \rightarrow S^1$, whose homotopy class will be
denoted by $n_i \in \pi_1(S^1)=\mathbb{Z}$ and will be referred to as the {\bf
local index} of $\overline{f}$ along $\Delta_i'$. It is easy to see that $n_i$
does not depend on the tubular neighborhood and the point $\delta_i$ used to
define it (but the sign of $n_i$ depends on orientations of the normal bundle
and $\Delta_i'$). If the lifting comes from a transversal section as above,
we certainly have $n_i=\pm 1$. Finally, we have the pullback $f^*\eta$ of
the canonical bundle over $S^2$, its Chern class will be denoted by 
$c_1 \in H^2(M \backslash \Delta)=H^2(M)$.

We now look at the map $\overline{f}$ as an ordinary defect map from $M$ to
$S^3$ with regular defect set $\Delta \cup \Delta'$. It is possible to
replace the surface $\Delta'$ by point defects by just altering the map
$\overline{f}$ on certain neighborhoods of the $\Delta_i'$: For each $i$,
choose a tubular neighborhood $N_i \subset M \backslash \Delta$ for 
$\Delta_i'$ with
$N_i \cap N_j= \emptyset$ for $i \neq j$, and let $DN_i$ 
denote the disk bundle,
$SN_i$ its boundary. Because we have $\pi_1(S^3)=\pi_2(S^3)=0$, we can find a
finite set $\Delta_i''$ for each $i$, lying in the interior 
$\smash {\overset{{}_\circ}{D}}N_i$ of $DN_i$, and a
defect map $\overline{g} : M \backslash (\Delta \cup \bigcup _i \Delta_i'')
\rightarrow S^3$ that coincides with $\overline{f}$ outside 
$\smash{\overset{{}_\circ}{D}}N = \bigcup_i \smash{\overset{{}_\circ}{D}}N_i$.
For each defect point $q \in \Delta'':=\bigcup_i \Delta_i''$ we now have
a local index $\iota_q(\overline{g}) \in \pi_3(S^3)=\mathbb{Z}$. 

In this example, the main question
that will be examined in this paper 
appears as the following problem: Are there any relations
between the points in $\Delta$ with their indices, the points in $\Delta''$
and some geometric data of the surfaces $\Delta_i'$? 
In this simple case, there is in fact a concrete formula, that will be
derived from a more general result presented later on:

\bigskip

\begin{satz} For each $i$, let $\chi_i$ be the self-intersection number
of $\Delta_i'$. Then 
$\sum_{q \in \Delta_i''} \iota_q(\overline{g})
= \pm n_i^2 \chi_i$, where the sign is the same for all $i$, and 
$\sum_i n_i^2\chi_i = c_1 ^2$.
\end{satz}

\bigskip

The sign in this formula does not depend on $M$, $\Delta'$ or $f$, it
only depends on the sign in a certain equation in the group $H^4(\Z,2;\Z)$ 
(see the proof for details).

Note that this result also applies to defect sections in $S^1$-bundles over
$M$, because up to point defects, every complex line bundle over $M$ can
be obtained as a pullback of the canonical bundle over $S^2$, associated to
the Hopf fibration.

As already mentioned above, Proposition 1 will be obtained as a special 
case of the main result of this paper that applies to the following 
situation:
Let $F \hookrightarrow E \stackrel{p}{\rightarrow} B$ be a fibration whose
fibre $F$ is of the type
$(\Pi,n)$ for some $n \geq 1$, 
i.e. $\pi_i(F)=0$ for $i \neq n$ and $\pi_n(F)=\Pi$, 
where $\Pi$
is an abelian group, and whose base space $B$ is a CW-complex.  Let $q$ be the
smallest integer with $\pi_q(E) \neq
0$, and suppose $q \geq n+2$. Using the homotopy sequence of the fibration, we
obtain $\pi_{n+1}(B)=\Pi$, the isomorphism given by the boundary map,
$\pi_i(B)=0$ for $1 \leq i \leq n$ and $n+2 \leq i \leq q-1$ and
$\pi_q(B)=\pi_q(E)$. Let $G$ denote the groups $\pi_q(E)$ and $\pi_q(B)$ that
will be identified in the sequel.

Now let $X$ be a connected CW-complex with a non-empty subcomplex $A$,
considered as
a relative CW-complex $(X,A)$ (using the terminology of Whiteheads book (\cite{W}),
so $X_k$ will denote the union of all $k$-cells not lying in $A:=X_{-1}$ with $A$).
Suppose there are given two continuous maps $f,g : X_q \rightarrow B$ that
coincide on $A$ and a partial lifting $k : A \rightarrow E$. 
\label{not} Now we have a well
defined primary obstruction $\gammaquer^{n+1}(k;f) \in H^{n+1}(X,A;\Pi)$ to
extend $k$ to a lifting of $f$ (originally, the obstruction is lying in the
group $H^{n+1}(X_q,A;\Pi)$, but will be considered as an element of the
isomorphic group $H^{n+1}(X,A;\Pi)$) and the corresponding obstruction
$\gammaquer^{n+1}(k;g)$ to extend $k$ as a lifting of $g$ (the coefficient
group really is $\Pi=\pi_n(F)$ because the base space $B$ is simply connected).
On the other hand, we have ordinary obstructions $\cquer^{q+1}(f),
\cquer^{q+1}(g) \in H^{q+1}(X,A;G)$ to extend $f$ and $g$ to the
$(q+1)$-skeleton. Finally, there is a unique cohomology operation $\Theta :
H^{n+1}( \cdot; \Pi) \rightarrow H^{q+1}(\cdot;G)$ with
$\Theta(\iota^{n+1})=k^{q+1}(B)$, where $\iota^{n+1} \in H^{n+1}(\Pi,n+1;\Pi)$
denotes the characteristic class and $k^{q+1}(B) \in H^{q+1}(\Pi,n+1;G)$ the
first non-trivial Postnikov invariant of $B$. We now have the following

\bigskip

\noindent {\bf Theorem} 
$\cquer^{q+1}(f) - \cquer^{q+1}(g) = \Theta(-\gammaquer^{n+1}(k;f)) - \Theta
(-\gammaquer^{n+1}(k;g)).$

\bigskip

There is a certain relationship between this equation and known results about
secondary obstructions. Consider for example the case $A=*$ and a map $f :
X_{i+1}
\rightarrow S^i$ defined on the $(i+1)$-skeleton of $X$, $i \geq 3$. 
By using an
embedding of $S^i$ into $K(\Z,i)$ (that may be constructed by attaching cells) 
and restricting the path-space fibration over $K(\Z,i)$ to $S^i$, 
we obtain
a fibration $E \rightarrow S^i$ as described above. If we choose the 
lifting $k$ to be the constant map and use the simple fact that $k^{i+2}(S^i)$ is 
the additive operation $Sq^2 \circ \rho_*$, where $\rho_*$ denotes 
reduction mod 2 (see the proof of
Proposition 2), we obtain
$$
\cquer^{i+2}(f)-\cquer^{i+2}(g)=Sq^2(\gammaquer^i(f)-\gammaquer^i(g) \text { mod
} 2).
$$
Using the definitions, it is easy to verify that the difference
$\gammaquer^i(f)-\gammaquer^i(g)$ is nothing else but the primary difference
$\dquer^i(f,g) \in H^i(X;\Z)$ (here again we make use of the boundary map in
the homotopy sequence to identify the coefficient groups), and we obtain 
the result
$$
\cquer^{i+2}(f)-\cquer^{i+2}(g)=Sq^2(\dquer^i(f,g) \text { mod
} 2)
$$
that has been proved by Steenrod in his paper \cite{St1}. 
So well known facts about
the secondary obstruction can be obtained as a special case of the result 
above. 
It should be possible to prove the Theorem by using known generalisations 
of Steenrods
result and the above relationship between the primary difference and the
lifting obstructions, but our proof avoids the difficulties that arise 
when the operation $\Theta$ fails to be additive, because
we do not make use of the primary difference between the maps $f$ and
$g$ into the base $B$.

\bigskip

Finally, a second example is considered, where the fibre is not of
Eilenberg-Mac Lane type, but the preceeding results can be applied with the 
help of an appropriate
decomposition of the fibration: Let $E \rightarrow S^4$ denote the unique
$SU(3)$-principal bundle whose second Chern class $c_2(P)$ is the orientation.
This bundle can be obtained by glueing together the trivial bundles over the
northern and southern hemisphere, using a representative of $1 \in
\pi_3(SU(3))=\mathbb{Z}$ as glueing map. It is well-known that
$\pi_1(SU(3))=\pi_2(SU(3))=0, \, \pi_3(SU(3))=\pi_5(SU(3))=\mathbb{Z}$ and
$\pi_4(SU(3))=0$ (see \cite{B}). From this, it is easy to deduce 
$\pi_i(E)=0$ for $i
\leq 4$ and $\pi_5(E)=H_5(E)=\mathbb{Z}$, using the exact homotopy sequence of
$E \rightarrow S^4$, the Wang sequence and the Hurewicz Theorem.

Now let $M$ be a manifold of dimension 6 and $f : M \rightarrow S^4$ a
continous map. Suppose there is a lifting $\fquer : M \backslash \Delta
\rightarrow E$ outside a closed surface $\Delta$. Again, we have a local index
for each connected component $\Delta_i$, 
but since we do not longer assume $M$ and $\Delta$ to
be orientable, the local index is well-defined only up to sign. Because of
that, we only consider its reduction mod 2, that will be denoted by $n_i$.
Because the total space $E$ is 4-connected, we again can replace the defect
set $\Delta$ by point defects if we just look at $\fquer$ as an ordinary map
into the total space $E$. Similar to the above, we can find a map 
$\gquer : M \backslash
\Delta' \rightarrow E$, where $\Delta'$ is a finite union of finite sets
$\Delta_i'$, each lying in the interior of a tubular neighborhood 
$DN_i \subset
M$ of $\Delta_i$, and $\gquer$ is supposed to coincide with $\fquer$ outside 
$\smash{\overset{{}_\circ}{D}}N = \bigcup_i \smash{\overset {{}_\circ}{D}}N_i$.
Finally we have local indices $\iota_q(\gquer)$ for each defect point $q \in
\Delta'$, and because
$\pi_5(E)=\mathbb{Z}$ and $M$ is not assumed to be oriented, this index is
well-defined in the group $\mathbb{Z}_2$.

\bigskip

\begin{satz} For each $i$, let $w_2(N_i) \in
H^2(\Delta_i,\mathbb{Z}_2)=\mathbb{Z}_2$ denote the second
Stiefel-Whitney-class of the normal bundle of $\Delta_i$. Then, for all $i$
$$
\sum _{q \in \Delta_i'} \iota_q(\gquer) = n_i \, w_2(N_i) \, \, \in
\mathbb{Z}_2.
$$
\end{satz}

\bigskip

\section {Proof of the Theorem}

\bigskip

During this proof, we again make use of the notations introduced on
page \pageref{not}.
Consider a Postnikov decomposition $\{ B^k,f_k \}$ of $B$, i.e. 
the space $B^k$ is obtained by attaching cells of dimension $\geq k+2$ at $B$,
$\pi_i(B^k)=0$ for $i >k$ and
there is a diagram
\begin{center} 
\begin{picture}(80,20) \usebox{\diagrammeins} \end{picture}
\end{center}
where the inclusions $\iota_k$ induce isomorphisms $\pi_i(B^k)=\pi_i(B)$ for $i
\leq k$. We can assume that the mappings $f_k$ are fibrations, whose fibres
$F_k$ are of the type $(\pi_k(B),k)$ and the inclusions $\pi_k(F_k) \rightarrow
\pi_k(B^k)$ are isomorphisms (see the discussion in \cite{W} p. 431). 
Furthermore it is
well known that $k^{q+1}(B)$ is just the primary obstruction to find a section
in $B^q \rightarrow B^{q-1}=K(\Pi,n+1)$.

Now consider the path space fibration $P=\{ \alpha : 
[0,1] \rightarrow B^{q-1} | 
\alpha(0)=f(*)=g(*) \} \rightarrow B^{q-1}$, where $* \in A$ is a
distinguished base point and the projection is given by sending a path
$\alpha$ to $\alpha(1)$ (see \cite{W} p. 31). It is clear that $P$ is 
contractible.
Using standard arguments (CW-approximation), we can assume $E$ to be a
CW-complex and apply obstruction theory to prove the existence of a lifting of
$\iota_{q-1} \circ p$, so we have the diagram
\begin{center}
\begin{picture}(80,30) \usebox{\diagrammzwei} \end{picture}
\end{center}
The  homotopy sequence  of the fibration 
$P \rightarrow B^{q-1}$ shows that the fibre $F'$ of 
$p'$
is of the type $(\Pi,n)$, where the isomorphism $\Pi  \cong \pi_n(F')$ is
induced by the restriction of $\eta$ to the fibre. An easy computation using
the definition of the boundary map
in the homotopy sequence proves that the primary obstruction to a section in $P
\rightarrow B^{q-1}$ is just $-\iota^{n+1} \in H^{n+1}(B^{q-1},\Pi)$.
Since $P$ is contractible, there is a homotopy
$$
\hdach : A \times [0,1] \rightarrow P
$$
with $\hdach_0=\eta \circ k$ and $\hdach_1=*$. Let $h:= p' \circ \hdach$,
so $h_0=\iota_{q-1} \circ f|A=\iota_{q-1} \circ g|A$ and $h_1=*$. 
Finally, let 
$\hquer : A \times [0,1] \rightarrow B^q$ be a lifting of $h$ with
$\hquer_0=\iota_q \circ f|A$. Since $\hquer_1$ and $*$ are two liftings of the
trivial map $*$, there is a well-defined primary difference
\begin{align}\label{d}
d:=\dquer^q(\hquer_1,*) \in H^q(A;G)
\end{align}
between these two maps as liftings
of $*:A \rightarrow B^{q-1}$.

Since
$B^{q-1}$ is an Eilenberg-Mac Lane space, we can extend $\iota_{q-1} \circ f$
over all of $X$, and using the homotopy extension property of $(X,A)$, we can
find a homotopy $H : X \times [0,1] \rightarrow B^{q-1}$ with $H_0 |
{X_q}=\iota_{q-1} \circ f$ that coincides with $h$ on $A \times [0,1]$. 
Let $f':=H_1 : (X,A)
\rightarrow (B^{q-1},*)$. 
Now we have a well-defined primary difference between
the two liftings $\hquer_1$ and $*$ of $f'|A=*$. 
Note that this difference is exactly the class $d \in H^q(A;G)$ defined in
equation (\ref{d}), since it only depends on the restriction
of $f'$ to $A$ which is the constant map $*$. 

There are primary obstructions $\gammaquer^{q+1}(\hquer_1;f')$ 
(resp. $\gammaquer^{q+1}(\hquer_0,H_0)$)to extend the
partial lifting of $f'$ ($H_0$), given by $\hquer_1$ ($\hquer_0$). 
From the coboundary formula (see \cite{St2}, 36.7) we have
$$
\gammaquer^{q+1}(\hquer_1;f')-\gammaquer^{q+1}(*;f')=\delta d,
$$
where $\delta : H^q(A;G) \rightarrow H^{q+1}(X,A;G)$ denotes the boundary.
Using this,
the homotopy invariance and the naturality properties of the 
obstructions and we can conclude 
\begin{align*}
\gammaquer^{q+1}(\hquer_0;H_0)=\gammaquer^{q+1}(\hquer_1;f') =&
\gammaquer^{q+1}(*;f') + \delta \dquer^q(\hquer_1,*) \\ 
=& \gammaquer^{q+1}(*;f') + \delta d= {f'}^*\gammaquer^{q+1}(*;id) + \delta d.
\end{align*}
But as mentioned above, the primary obstruction $\gammaquer^{q+1}(*;id)$ is
nothing else but
the Postnikov invariant, and it is well known (see \cite{W} p. 450) that the class
$\gammaquer^{q+1}(\hquer_0;H_0)=\gammaquer^{q+1}(\iota_q \circ f; H_0)$
coincides with the obstruction $\cquer^{q+1}(f)$ (remember that $H_0$ is
an extension of $\iota_{q-1} \circ f$). So we obtain the relation
$$
\cquer^{q+1}(f) = {f'}^* k^{q+1}(B) + \delta d.
$$
On the other hand, we can consider the fibration $P \rightarrow B^{q-1}$ with
primary lifting obstruction $-\iota^{n+1} \in H^{n+1}(X,A;\Pi)$, and we have
$$
{f'}^* \iota^{n+1} = - {f'}^* \gammaquer^{n+1}(*;id) = -
\gammaquer^{n+1}(*;H_1) = - \gammaquer^{n+1}(\eta \circ k ; H_0),
$$
but the latter class is (remember that we identified the homotopy of the
fibre of $P \rightarrow B^{q-1}$ with that of $F$, using the map $\eta$)
exactly $-\gammaquer^{n+1}(k;f)$, so we finally have
$$
{f'}^*\iota^{n+1} = - \gammaquer^{n+1}(k;f).
$$
Together with the above relation, we obtain
$$
\Theta(-\gammaquer^{n+1}(k;f))= \Theta({f'}^*\iota^{n+1})={f'}^*k^{q+1}(B) =
\cquer^{q +1}(f) - \delta d.
$$
Replacing $f$ by $g$ in the above construction yields the
formula 
$$
\Theta(-\gammaquer^{n+1}(k;g))= \cquer^{q +1}(g) - \delta d
$$
(note that we can use the same homotopies $\hquer$ and $\hdach$ and so
we have the same class $d$, since $f$ and $g$ coincide on $A$).
Subtraction of these equations give the desired result. \hfill $\square$

\bigskip

\section {Proof of Proposition 1}

\bigskip

Now we specialise to the case where $E \rightarrow B$ is the Hopf fibration,
$q= 3,n=1$ and the pair $(X,A)$ is $(DN_i,SN_i)$, with the notations already
used in the
introduction ($i$ will be fixed in the sequel). The map $f$ is the
restriction of the point-defect map $f :
M \backslash \Delta \rightarrow S^2$. Since $DN_i$ does not contain defect
points, we have $\cquer^{q+1}(f)=0$. Let $g:=\pi \circ \gquer$, where $\pi$
denotes the projection $S^3 \rightarrow S^2$. It is easy to see (using the
identification $\pi_3(S^3)=\pi_3(S^2)=\mathbb{Z}$ induced by $\pi$ and a
triangulation of $DN_i$ whose 3-skeleton does not intersect $\Delta_i''$), that
the obstruction $\cquer^{q+1}(g) \in H^4(DN_i,SN_i)=\mathbb{Z}$ is just the 
sum  $\sum_{q \in \Delta_i''} \iota_q(\gquer)$. 
Now choose orientations of the
normal bundle and $\Delta_i'$ with the property that under the
Thom isomorphism, the orientation of the surface $\Delta_i'$ corresponds to the
orientation of $DN_i$, induced by that of $M$ (since the formula contains the 
square $n_i^2$, we
can suppose that these orientations were used to define the index $n_i$). 
Then the class $\gammaquer^{n+1}(k;f)$ is a certain multiple of the 
Thom class $\tau_i$, in fact it is $n_i \tau_i$, and certainly 
$\cquer^{n+1}(k;g)=0$. 
Finally, the operation $\Theta$ associated to the Postnikov invariant 
$k^{q+1}(S^2)$ is just the cup square, up to sign. Since the square 
$\tau^2_i$ is $\chi_i [DN]$, where $[DN]$ denotes the orientation, 
we obtain the first part of Proposition 1. The second part is
just the fact that certainly $\sum_{p \in \Delta} \iota_p(f)
+ \sum_{q \in \Delta'} \iota_q(\gquer)=0$, and the first sum
is just $\pm c_1^2$, where the sign is the same as before (this follows
if we apply the Theorem with $(X,A)=(M,*)$ and $g=*$).
\hfill $\square$

\bigskip

\section {Proof of Proposition 2}

First, we have to compute the Postnikov invariant $k^6(S^4) \in
H^6(\Z,4;\Z_2)$. It is a result of Eilenberg and Mac Lane (\cite{EM} p. 122) 
that
this group is $\Z_2$. The non-zero element is given by the operation $Sq^2
\circ \rho_*$, where $\rho$ denotes the reduction mod 2: Because the operation
is stable, $Sq^2 \circ \rho_* = 0$ in dimension 4 would imply $Sq^2 \circ
\rho_* = 0 : H^2(X;\Z) \rightarrow H^4(X;\Z)$, what is easily seen not to be
true, just consider the example $X=\mathbb{C} P^{\infty}$.  A similar argument,
using a suspension formula for the Postnikov invariants (\cite{W}
p. 439), tells us that the operation $k^6(S^4)$ is not zero (this
would give $\rho_*k^4(S^2)=0$, but $k^4(S^2)$ is just the cup square). 
Hence $k^6(S^4)=Sq^2 \circ \rho$
(in fact we have proved $k^{i+2}(S^i)=Sq^2 \circ \rho_*$ for all $i \geq 4$).

We now decompose our fibration using the following Lemma
(that can be seen as the first stage of the so called
``Moore-Postnikov decomposition'' ,see \cite{S}, Chapter 8 for details):

\bigskip

\noindent
{\bf Lemma 1:} Let  $F \hookrightarrow E \stackrel{p}{\rightarrow} B$
be a fibration with $(n-1)$-connected $E$, $n \geq 2$, simply connected
base space $B$ and simple and connected fibre $F$. The total space $E$ 
is supposed to be of
the homotopy type of a CW-complex. Let $q$ be an integer , $1\leq
q <n$.

Then there are fibrations  $F' \hookrightarrow S
\stackrel{\pi}{\rightarrow} B$
and $F'' \hookrightarrow \tilde{E} \stackrel{q'}{\rightarrow} S$, together with
maps $q : E \rightarrow S$ and $\lambda : E \rightarrow \tilde{E}$
such that the diagram   
\begin{center}
\begin{picture}(30,40) \usebox{\diagrammdrei} \end{picture}
\end{center}
commutes and the following holds:
\begin{enumerate}
\item $\lambda$ is a homotopy equivalence,
\item $q_* :  \pi_i(F) \rightarrow \pi_i(F')$ is an 
isomorphism for $i \leq q-1$ and  $\pi_i(F')=0$ for $i \geq q$,
\item $S$ is $q$-connected, 
$\pi_* : \pi_i(S) \rightarrow \pi_i(B)$ is an isomorphism 
for  $i \geq q+1$ , and
\item  $\pi_i(F'')=0$ for $i \leq q-1$, and 
$\pi_i(F'') \cong \pi_i(F)$ for $i \geq q$ .
\end{enumerate}

\bigskip

\noindent
{\bf Proof:} By attaching cells of dimension $\ge q+2$, we can construct a
cellular extension $B^*$ of $B$ with $\pi_i(B^*)=0$ for $i \geq q+1$. Let $P'
\rightarrow B^*$ be the path space fibration and $S$ be the total space of its
restriction to $B$. It is easy to see that there is a lifting $q : E \rightarrow
S$ of $E \rightarrow B$, using the assumptions made on the homotopy groups of
$E$, and the existence of the claimed homotopy equivalence $\lambda$ and the
fibration $\tilde{E} \rightarrow S$ is based on a standard tool in homotopy
theory, see \cite{W} p. 42. It is now a straightforward computation
to verify the above
properties, using the homotopy sequences for the involved fibrations. \hfill
$\square$

\bigskip

Now choose a decomposition as in the Lemma for the $SU(3)$-bundle $E
\rightarrow S^4$ with $q=4$. Then $\pi_5(S)=\pi_5(S^4)=\Z_2$, the fibre $F'$
is of the type $(\Z,3)$, and with
the identifications $\pi_5(E)=\Z$, $\pi_5(S)=\Z_2$, the map $q_* : \pi_5(E)
\rightarrow \pi_5(S)$ is easily seen to be an epimorphism, so it is just the
canonical map $\Z \rightarrow \Z_2$. We can suppose that $\Delta$ is connected
and  $M=DN$. Now let $\gammaquer^4(q \circ \fquer;f) \in H^4(DN,SN;\Z)$ be the
primary obstruction to extend the lifting $q \circ \fquer$ of $f$ to $DN$. It is
easy to see that $\gammaquer^4(q \circ \fquer;f) \text{ mod } 2 = n \tau$,
where $\tau$ denotes the $\Z_2$ Thom class of the normal bundle $N$. Let
$\cquer^6(q \circ \fquer) \in H^6(DN,SN;\Z_2)$ denote the primary obstruction
to extend $q \circ \fquer$ as an ordinary map into $S$. Applying the above
calculation of $k^6(S^4)$ and the Theorem, similar to the proof of
Proposition 1, we obtain $\cquer^6(q \circ \fquer) = Sq^2(n\tau)=n Sq^2(\tau)$.
But clearly $\cquer^6(q \circ \fquer) \in H^6(DN,SN;\Z_2)=\Z_2$ is just the sum
$\sum_{q \in \Delta'} \iota_q(\gquer)$, and because $Sq^2(\tau)$ is identified
with the second
Stiefel-Whitney class $w_2(N)$ under the Thom isomorphism (see \cite{H}, 
Chapter 17, 9.1),  we obtain the desired result
$\sum_{q \in \Delta'} \iota_q(\gquer)=n w_2(N)$. \hfill $\square$

\bigskip

The author  
would like to thank Klaus J\"anich and Stefan Bechtluft-Sachs
for their support and helpful discussions during the work on this paper.

\end{document}